\documentclass[12pt]{amsart}

\usepackage{a4wide}

\usepackage{amsmath, amsthm, amsfonts, amssymb,accents,a4}
\usepackage{dsfont}
\usepackage{srcltx, color}

\theoremstyle{plain}
\newtheorem{theorem}{Theorem}[section]
\newtheorem{lemma}[theorem]{Lemma}
\newtheorem{corollary}[theorem]{Corollary}
\theoremstyle{remark}
\newtheorem{remark}[theorem]{Remark}

\numberwithin{equation}{section}

\def\om{\omega}
\def\e{\varepsilon}
\def\g{\gamma}
\def\G{\Gamma}
\def\l{\lambda}
\def\p{\partial}
\def\D{\Delta}

\def\E{\mbox{\rm e}}

\def\d{\delta}

\def\H{W_2}

\def\di{\,\mathrm{d}}

\def\I{\mathrm{I}}

\def\la{\langle}
\def\ra{\rangle}

     \newcommand{\CC}{\mathds{C}}

     \newcommand{\EE}{\mathbb{E}}

     \newcommand{\NN}{\mathds{N}}

     \newcommand{\PP}{\mathbb{P}}

     \newcommand{\RR}{\mathds{R}}
     \newcommand{\ZZ}{\mathds{Z}}

\DeclareMathOperator{\dist}{dist}

\begin{document}

\allowdisplaybreaks

\title[Spectra of weak-disorder quantum waveguides]{Low lying spectrum of weak-disorder quantum waveguides}
\author[Borisov and Veseli\'c]{Denis Borisov$^{1,2,3}$ and  Ivan Veseli\'c$^{1,4}$}
\maketitle

\begin{quote}
\begin{itemize}
\item[1)] \emph{Faculty of Mathematics, Technische Universit\"at Chemnitz, 09107 Chemnitz, Germany}
\item[2)] \emph{Department of Physics and Mathematics, Bashkir
State Pedagogical University, October rev. st.~3a, Ufa, 450000,
Russia}
\item[3)] \emph{E-mail:}\ \texttt{borisovdi@yandex.ru}, \emph{URL:} \texttt{http://borisovdi.narod.ru/}
\item[4)] \emph{URL:}\ \texttt{http://www.tu-chemnitz.de/mathematik/stochastik/}
\end{itemize}
\end{quote}

\begin{abstract}
We study the low-lying spectrum of the Dirichlet Laplace operator on a
randomly 
 wiggled strip. More precisely, our results are formulated in terms of the eigenvalues of finite segment approximations of the infinite waveguide.
Under appropriate weak-disorder assumptions we obtain deterministic and probabilistic bounds on the position of the lowest eigenvalue.
A Combes-Thomas argument allows us to obtain a so-called 'initial length scale decay estimates' at they are %�%used
employed in the proof
of spectral localization using the multiscale analysis method.
\end{abstract}

\section{Introduction}
The propagation of waves in disordered  media can be modeled using differential equations governed by a random Hamiltonian. The most important questions in this context concern the long time behavior of
waves which oftentimes 
allows to conclude results concerning the transport properties of the material described by the random operator. A particularly well studied class of operators
is the one which arises in the quantum mechanical description of disordered solids.
To this class belong various types of random Schr\"odinger operators, e.g.{} the Anderson model
on $\ell^2(\ZZ^d)$, or a Laplacian with Poisson distributed repulsive impurity potentials on $L^2(\RR^d)$.

For most operators it is possible to relate propagation properties
to spectral features, e.g.{} by the use of the RAGE or Ruelle theorem. From this point of view it is justified to study first the measure theoretic
spectral types which arise in a certain random model, and then relate them
to the transport properties of the considered material.
This leads to the study of a plethora of spectral features, some of which are specific to the situation  that we are dealing with a random operator, i.e.~with an infinite family of individual operators.
Let us mention some of these properties: the characteristics of the spectrum as a subset of the real line, its band structure, 
location of the spectral infimum (i.e.{} the infimum of the
bottom of the spectrum when varied over all members of the family),
the density of the spectrum in various energy regions, decay properties of the Green's function, etc.

To be able to point out the interesting contribution of the paper at hand, we
shall assume in the further discussion
that the randomness enters the Hamiltionan via
a countable family of random variables.
Already when considering the very basic features of the spectrum,
one sees that it makes a great difference whether the operator (more precisely: the associated quadratic form) depends in a monotone or non-monotone way on the random variables entering the model.
In the case of monotone dependence the spectral minimum of the operator family obviously corresponds
to the configuration where all  random variables  are set to one of the extremal values. Similarly, in the monotone situation the band structure of the spectrum can be analyzed 	
using rather basic sandwiching arguments,
see e.g.{} \cite{KirschSS-98a}. It is consistent with these elementary examples of the advantages of monotonicity
that there is a rather good understanding of typical energy/disorder regimes where  monotone models exhibit
localization of waves, see the monographs and survey articles
\cite{Stollmann-01,KirschM-07,Veselic-07b,Kirsch-08}.

If the dependence of the quadratic form on the  random variables  is not monotone,
already the identification of the spectral minimum is not obvious and sometimes a highly non-trivial question, see e.g.{}
\cite{BakerLS-08,KloppN-09a}.
For more intricate properties, like the regularity of the density of states
or the analysis of spectral fluctuation boundaries, the difference between monotone
and non-monotone models is even more striking.

Nevertheless there has been a continuous effort to bring the understanding of models
with non-monotone dependence on the randomness to a similar level as the one for monotone models.
The model of this type to which most attention was devoted so far is the alloy type Hamiltonian with single site potentials of changing sign,
see e.g.{} \cite{Klopp-95a,Stolz-00,Klopp-02c,Veselic-02a,HislopK-02,KostrykinV-06, KloppN-09a}.
More recently also the discrete analog of this model was studied in \cite{Veselic-circulant,ElgartTV,Veselic,TautenhahnV-10,Krueger}.
Electromagnetic Schr\"odinger operators with random magnetic field
 \cite{Ueki-94,Ueki-00,HislopK-02,KloppNNN-03,Ueki-08, Bourgain-09}, as well as Laplace-Beltrami
operators with random metrics \cite{LenzPV-04,LenzPPV-08,LenzPPV-09} are other examples with non-monotonous parameter dependence.

A very interesting model with geometric disorder is the random displacement model,
which exhibits also no obvious monotonicity, cf.{} e.g.{}
\cite{Klopp-93,BakerLS-08,KloppLNS}.
Another relevant model (although not defined in terms of a countable family of i.i.d.{}  random variables)
without obvious monotonicity is a random potential given by a Gaussian stochastic field with sign-changing covariance function, c.f.{} \cite{HupferLMW-01a,Ueki-04,Veselic-gauss}.

In this paper we consider a family of Hamiltonians which consists of  the Dirichlet Laplacian on
a randomly  
wiggled waveguide. In this model the dependence of the quadratic form on the  random variables
is neither monotone nor linear. In this respect it is related to the random displacement model.
Moreover, in our model, the randomness does not enter via potential terms, but rather
through differential operator terms.
A random waveguide model has been studied before in \cite{KleespiesS-00}.
There the randomness enters via a variation of the width of the waveguide.
This type of perturbation leads to a quadratic form which depends monotonously on the  random variables  and is thus structurally different from our model.

There is a substantial body of literature devoted to the analysis
of eigenvalues below the essential spectrum of bent, asymptotically straight waveguides, see for instance \cite{ExnerS-89}.
These eigenvalues are in contrast to the purely absolutely continuous spectrum exhibited by an straight waveguide.
Thus a local geometric perturbation leads to the emergence of discrete eigenvalues.
Given this fact, it is interesting to ask whether geometric perturbations which are ergodic and random
lead to 
dense point spectrum below the continuous one,  in analogy to the phenomenon encountered for several classes of random  Schr\"odinger operators mentioned above.
One should point out that in the present paper the local geometric perturbations of the waveguide
are introduced in a somewhat different way than in \cite{ExnerS-89}.

Let us now describe the main result of this paper.
We derive lower bounds on the first eigenvalue for a finite segment of a randomly 
wiggled strip in $\RR^2$.
They measure how far  the eigenvalue may move up, if the vector of random variables moves away from the optimal configuration.
As an application we obtain a second result:
In the terminology of the \emph{multiscale analysis} (MSA)
it is a \emph{initial length scale estimate}
for energies near the bottom of the spectrum in the weak disorder regime.
It corresponds to the induction anchor of the
MSA. This should be understood as a step towards a localization proof using the MSA.
If there would be an appropriate Wegner estimate at disposal at low energies (which we don't have
at the moment) an adaptation of the usual MSA, e.g. as presented in \cite{Stollmann-01,GerminetK-01a,GerminetK-04}
would lead to localization.

Let us say a few words about our methods of proof. It consists of a deterministic
and probabilistic part. For a finite segment of the waveguide one can use methods
from asymptotic analysis to estimate the position of the principal eigenvalue of the Laplacian.
In this situation only a finite number of  random variables  enters the operator.
It is this part which requires the weak disorder restriction.
The mentioned results can be combined
with a Dirichlet-Neumann bracketing argument and
a large deviations principle to arrive at an exponential
probabilistic bound on the position of the lowest finite segment eigenvalue.
Using a Combes-Thomas estimate \cite{CombesT-73,BarbarouxCH-97b, Stollmann-01,BoutetdeMonvelS-03b}
this can be turned into an off-diagonal decay estimate on the Green's function,
which plays the role of the initial length scale estimate in the MSA.

It is maybe worthwhile to point out some differences to the recent paper \cite{KloppN-09a} of Klopp and Nakamura
which is devoted to the proof of Lifshitz tails for alloy-type Schr\"odinger operators with
single site potentials which are allowed to change sign.
There are two aspect in common between this work and ours: both of them
concern the analysis of the low lying eigenvalues of finite volume random Hamiltonians
and both of them deal with non-monotone parameter dependence. There are also two main differences:
we are not able to give an Lifshitz  bound on the integrated density of states for our model, since
the global disorder coupling constant has to be 
chosen dependent on the volume scale. If one lets the scale tend to infinity the
global coupling constant has to go to zero. On the other hand we assume no reflection symmetry for the individual perturbations.
A crucial assumption of  \cite{KloppN-09a}  is that the single site potential obeys such an condition. It allows Klopp and Nakamura to use,
after some work and  ingenious ideas,
an effective  decoupling between different random variables (similarly as in the case of fixed-sign single site potentials).
This means that it is only necessary to perform perturbation theory with respect to one coupling constant.
This aspect of the proof of  \cite{KloppN-09a}  is discussed on page 1134 before the statement of hypothesis (H2) there.
In our model there is no such symmetry assumption. This means that the analysis of the single parameter random Hamiltonian
on a unit cell with Neumann b.c. does not give us the crucial information which was instrumental in the proof strategy of  \cite{KloppN-09a}.
Consequently,  we need to analyze the fully interacting model, which leads to an eigenvalue perturbation problem with respect to
many parameters.

The paper is organized as follows: in the next section we define rigorously our model,
in Section \ref{s:probab-shifted} we state the probabilistic estimates on the position of the principal
eigenvalue of a finite segment of a random waveguide and on the exponential off-diagonal decay of the associated Green's function,
in Section \ref{s:reduction} we reduce the proof of the two above statements to deterministic bounds on the first eigenvalue,
and in the final Section \ref{s:deterministic} the mentioned deterministic estimates are established.

\section{Model}
\label{s:model}

We consider random quantum waveguides in $\RR^2$ ,
determined by the following data:
Let $(\omega_k)_{k\in \ZZ}$ be a sequence of independent, identically distributed, non-negative, bounded, non-trivial random variables,
$\kappa>0$ a global coupling constant, $l\geqslant 1$ the length of one (periodicity) cell of the waveguide,
and $g \in C_0^2(0,l)$ a single bump function.
The following function $G\colon\RR \times \Omega \to \RR$ determines the shape of the waveguide
\[
 G(x_1,\omega) := \sum_{k\in\ZZ} \omega_k \, g(x_1 -kl ) .
\]
Note that $\kappa G(x_1,\omega) =G(x_1,\kappa \omega)$.
Together with the global coupling constant $\kappa >0$ it defines
an infinite  waveguide as the set
\begin{align*}
D_{\kappa,\omega}
:& =\{x \in \RR^2 \mid x_1\in\RR, \kappa G(x_1,\omega)<x_2<\kappa G(x_1,\omega)+\pi\}
\\
& =\{x \in \RR^2 \mid x_1\in\RR, G(x_1,\kappa \omega)<x_2< G(x_1,\kappa \omega)+\pi\}
=
D_{1,\kappa\omega}
\end{align*}
Our results  are  formulated in terms of spectral features of
finite segments of the infinite waveguide. We define them next.
For $N\in\NN$ and  $j\in\ZZ$ set
\[
D_{\kappa,\omega}(N,j):=\{x \in \RR^2 \mid jl< x_1< (j+N)l , \kappa G(x_1,\omega)<x_2<\kappa G(x_1,\omega)+\pi\}.
\]
Denote  by $\Gamma_{\kappa,\omega}(N,j)$ the upper and lower part of the boundary of
$D_{\kappa,\omega}(N,j)$, i.e.,
\begin{align*}
\Gamma_{\kappa,\omega}(N,j):=&
\{x \in \RR^2 \mid jl< x_1< (j+N)l , x_2=\kappa G(x_1,\omega)\}
\\ \cup &
\{x \in \RR^2 \mid jl< x_1< (j+N)l , x_2=\kappa G(x_1,\omega)+\pi\}.
\end{align*}
The remaining part of the boundary $\partial D_{\kappa,\omega}(N,j)\setminus
\Gamma_{\kappa,\omega}(N,j) $ is denoted by   $\gamma_{\kappa,\omega}(N,j)$.

Let
$\mathcal{H}_{\kappa,\omega}(N,j)$
denote
the negative Laplace operator on
$D_{\kappa,\omega}(N,j)$ with Dirichlet boundary conditions on $\Gamma_{\kappa,\omega}(N,j)$
and Neumann b.c.~on $\gamma_{\kappa,\omega}(N,j)$. The lowest eigenvalue of
$\mathcal{H}_{\kappa,\omega}(N,j)$ is denoted by $\lambda_{\kappa,\omega}(N,j)$.
If $j=0$ we shall 
 use the following shorthand notation:
\begin{align*}
D_{\kappa,\omega}(N)&:=D_{\kappa,\omega}(N,0),&
\Gamma_{\kappa,\omega}(N)&:=\Gamma_{\kappa,\omega}(N,0),&
\gamma_{\kappa,\omega}(N):=\gamma_{\kappa,\omega}(N,0),
\\
\mathcal{H}_{\kappa,\omega}(N)&:=\mathcal{H}_{\kappa,\omega}(N,0),&
\lambda_{\kappa,\omega}(N)&:=\lambda_{\kappa,\omega}(N,0).&
\end{align*}
Similarly as for the infinite waveguide we have
$D_{\kappa,\omega}(N,j) = D_{1,\kappa\omega}(N,j)$.
Since $\kappa >0$ is arbitrary we may assume without restricting the model
\begin{equation}\label{1.0b}
\max\{\|g\|_{C[0,l]}, \|\|g'\|_{C[0,l]}, \|g''\|_{C[0,l]}\}=1.
\end{equation}

Denote the distribution measure of $\omega_k $ by $\mu$.
It will be convenient to think of $\mu$ as a measure on the  semiaxis $[0,\infty)$
with support in the unit interval $[0,1]$.
Thus any $\omega_k$ takes values larger than $1$ only with zero probability.
Then $\PP = \bigotimes\limits_{k\in\ZZ}\mu$ denotes the
product measure  on the configuration space
$\Omega=\times_{k\in\ZZ}[0,\infty)$ whose elements we
denote by $\omega = (\omega_k)_{k \in\ZZ}$.

Note that by the assumptions on $\mu$ the following statements hold for $\PP$-almost all $\omega \in \Omega$: $\omega \in \ell^\infty(\ZZ)$
and $\|\omega \|_\infty := \sup_{k \in \ZZ} |\omega_k| =\sup_{k \in \ZZ} \omega_k <\infty$.
On appropriate subspaces of $\Omega$
we define the usual  $\ell^p$-norms, in particular
$\|\omega\|_1:= \sum\limits_{k \in \ZZ} \omega_k$
and $\|\omega\|_2:= (\sum\limits_{k \in \ZZ} \omega_k^2)^{\frac{1}{2}} $.

 The randomness of the finite segments $D_{\kappa,\omega}(N,j)$
arises from only a finite number of  random variables  $(\omega_k)$. For this reason it is
convenient to have the following notation at disposal: For $\Lambda \subset \ZZ$
define the projection map
\[
\pi_\Lambda \colon \Omega \to \Omega , 
\quad (\pi_\Lambda(\omega))_k =\omega_k \, \chi_\Lambda (k)
\]
Here $\chi_A$ denotes the indicator function of a set $A$.
If $\Lambda \subset \ZZ$ is finite  $\pi_\Lambda (\Omega)$ is contained $\ell^p(\ZZ)$
for  every  $p\in [1,\infty]$. We 
 shall use the shorthand notation
$\|\omega\|_{\Lambda,p} := \| \pi_\Lambda \omega\|_{p} $
and in the case $\Lambda= \{1,\dots,N\}$
\[
\|\omega\|_{N,p}  := \|\omega\|_{\Lambda,p} := \| \pi_\Lambda \omega\|_{p}.
\]

\section{Probability of low lying eigenvalues and the initial length scale estimate}
\label{s:probab-shifted}

Here we present estimates on the probability that the lowest eigenvalue of
$\mathcal{H}_{\kappa,\omega}(N,j)$ is close to one.
The first information which one has is that the minimum of the spectrum of
the Laplacian on an straight waveguide segment $\mathcal{H}_{0,\omega}(N,j)$ is equal to one.
This follows directly from separation of variables.
We 
shall see later that
no operator $\mathcal{H}_{\kappa,\omega}(N,j)$  has spectrum below one.
In this sense we can say that one is the minimal spectral value
for all random configurations.

To formulate the main result we 
shall use the following quantities:
\begin{equation*}
\widetilde{g}:=g-\frac{1}{l}\int_0^l g(t)\di t,\quad
 c_2 =\frac{3\,\|\widetilde{g}\|_{L_2(0,l)}}{2\, l^3},
 \quad c_3 =\frac{3\,\|\widetilde{g}\|_{L_2(0,l)}^2}{5000\, l^7}.
\end{equation*}

\begin{theorem}
\label{th:ilse}
Let $g$ and $\mu$ as above be given. Let $\gamma >34$.
Then there exists an initial scale $N_1$ such that if $N \geqslant N_1$  the interval
\[
 I_N:=\left[\frac{2 N^{\frac{1}{\gamma}-\frac{1}{4}}}{\EE\{\omega_k\}\sqrt{c_2}} , c_3  N^{-\frac{15}{2\gamma}} \right]
\]
is non-empty. If $N \geqslant N_1$  and $\kappa \in I_N$, then
\begin{equation}
\label{eq:initial}
 \PP\left(\omega \in \Omega\mid \lambda_{\kappa,\omega}(N)-1\leqslant N^{-\frac{1}{2}}\right )
\leqslant
N^{1-\frac{1}{\gamma}} \ \E^{-c_4N^{1/\gamma}}
\end{equation}
for a constant $c_4>0$ depending only on $\mu$.
\end{theorem}

\begin{remark}
The definition of the interval $I_N$ encodes  how  our  weak-disorder regime depends on the length scale $N$.
As an instance let us choose $\gamma =35$. Then it is possible to choose $\kappa = c_3 N^{-\frac{1}{4}}$,
which means that the allowed disorder regime does not shrink very fast for $ N \to \infty$.
For larger $\gamma$ the behavior is even better.
\end{remark}

Using a Combes-Thomas estimate we arrive at the following estimate
on the probability that the Green's function, resp. resolvent,  exhibits exponential off-diagonal decay
for energies very close the energy one, i.e.{} the overall minimum of the spectrum.
Note that for $N$ large $I_N \subset [0,1]$.

\begin{corollary}\label{c:Combes-Thomas}
Let $g$, $\mu$, $I_N$, $\gamma$, $N_1$ and $c_4$ be as in Theorem \ref{th:ilse}.
Let $\kappa \in I_N \cap [0,1]$, $\alpha, \beta \geqslant 2$ and set
\begin{align*}
A  &:=\{x \in \RR^2 \mid 0\leqslant  x_1\leqslant  \alpha  , \kappa G(x_1,\omega)<x_2<\kappa G(x_1,\omega)+\pi\} \subset D_{\kappa, \omega} (N),
\\
B  &:=\{x \in \RR^2 \mid L-\beta\leqslant  x_1\leqslant  L  , \kappa G(x_1,\omega)<x_2<\kappa G(x_1,\omega)+\pi\} \subset D_{\kappa, \omega} (N).
\end{align*}
Then we have for any $N \geqslant N_1$
\begin{multline*}
 \PP\left(\omega \in \Omega\mid \forall  \ \lambda\in [1,1+1/(2\sqrt N)]  \colon
\left \|\chi_A (\mathcal{H}_{\kappa,\omega}(N)-\lambda)^{-1}   \chi_B \right \|
\leqslant
2\sqrt N
\E^{\left({-\frac{\dist(A,B)}{48
\, \sqrt N}}\right)}
\right)
\\
\geqslant
1-N^{1-\frac{1}{\gamma}} \ \E^{-c_4N^{1/\gamma}}.
\end{multline*}
\end{corollary}

In the typical formulation of an initial scale estimate for the multiscale analysis
one requires the probability of
\[
\left \{\exists  \ \lambda\in [1,1+1/(2\sqrt N)]  \colon \left \|\chi_A (\mathcal{H}_{\kappa,\omega}(N)-\lambda)^{-1}   \chi_B \right \|
>
2\sqrt N  \E^{-\frac{\dist(A,B)}{48\,\sqrt N}} \right\}
\]
to be bounded by an inverse power of $N$.
Let $ q \in \NN$ and $ N_0 $ be such that $N^{(1+q-\frac{1}{\gamma})/c_4} \leqslant \exp(N^{\frac{1}{\gamma}})$ for all $N \geqslant  N_0 $.
Set $N_2:= \max (N_1, N_0 )$. Then for any $N \geqslant N_2$ we have
\begin{multline*}
\PP\left(\omega \in \Omega\mid \forall  \ \lambda\in [1,1+1/(2\sqrt N)]  \colon \left \|\chi_A (\mathcal{H}_{\kappa,\omega}(N)-\lambda)^{-1}   \chi_B \right \| \leqslant 2\sqrt N \E^{-\left({-\frac{\dist(A,B)}{ 48 \, \sqrt N}}\right)} \right) \\ \geqslant 1- N^{-q}.
\end{multline*}

\section{Proof of Theorem \ref{th:ilse} and Corollary  \ref{c:Combes-Thomas}}
\label{s:reduction}
In this section we prove that Theorem \ref{th:ilse}
is implied  by the following Theorem and its Corollary.
We also show how  Corollary~\ref{c:Combes-Thomas}
follows.

\begin{theorem}
\label{th:deterministic-lower-bound-1}
Recall that
$(\omega_k)_{k\in\ZZ}$ is a sequence of non-negative reals.
Let $g$, $l$, $\mathcal{H}_{\kappa,\omega}(N)$, $\lambda_{\kappa,\omega}(N)$ be as in Section \ref{s:model}.
Recall that
\begin{equation*}
\widetilde{g}:=g-\frac{1}{l}\int_0^l g(t)\di t.
\end{equation*}
Assume that
\begin{equation*}
\kappa \|\omega\|_{N,2}
\leqslant   \frac{3}{5000}  \, \|\widetilde{g}\|_{L_2(0,l)}^2  \, \frac{1}{ l^7N^7}= \frac{c_3}{N^7}
\end{equation*}
Then the estimate
\begin{equation*}
\lambda_{\kappa ,\omega}(N)-1
\geqslant
\frac{3}{2}\,  \|\widetilde{g}\|_{L_2(0,l)}^2  \, \frac{\kappa^2\|\omega\|_{N,2}^2}{l^3N^3}
=c_2 \frac{\kappa^2 \|\omega\|_{N,2}^2}{N^3}
\end{equation*}
holds true.
\end{theorem}

\begin{corollary}
\label{th:lower-bound}
If $\kappa    \leqslant   c_3 N^{-\frac{15}{2}}$, then we have for $\PP$-almost all $\omega \in \Omega$:
\[
\lambda_{\kappa,\omega}(N) -1 \geqslant
c_2 \frac{\kappa^2 \|\omega\|_{N,2}^2}{N^3}.
\]
\end{corollary}
\begin{proof}
 This follows immediately from Theorem \ref{th:deterministic-lower-bound-1}, since $\|\omega\|_{N,2} \leqslant \sqrt{N} \|\omega\|_{N,\infty} \leqslant \sqrt{N} $
for almost all $\omega$.
Thus $\kappa \leqslant constant\,N^{-\frac{15}{2}}$  implies $\kappa \|\omega\|_2 \leqslant constant\, N^{-7}$.
\end{proof}

The following \emph{large deviations principle} will be used in the proof of Theorem \ref{th:ilse}.

\begin{lemma}\label{l:LDP}
Let $\omega_k, k\in\ZZ$ be an i.i.d.~sequence of non-trivial, non-negative,
bounded random variables. Then there exists a constant $c_4>0$ depending only on $\mu$ such that
\[
 \forall \ n \in \NN \ : \
\PP\left(\omega \mid \frac{1}{n} \sum_{k=1}^n \omega_k \leqslant \frac{\EE\{\omega_k\}}{2} \right ) \leqslant \E^{-c_4n} .
\]
\end{lemma}

\begin{proof}[Proof of Theorem \ref{th:ilse}]
For $ K \in 2 \NN$ and $\gamma \in \NN$ we set $N:= K^\gamma$.
Note that we can decompose a long waveguide segment into smaller parts. Thus, up to a set of measure zero
$D_{\kappa,\omega}(N)  $ equals
\[
\bigcup\limits^\bullet_{j=0, \dots , J-1} D_{\kappa,\omega}(K,j)
\]
where $J =N/K = K^{\gamma-1}=N^{1-\frac{1}{\gamma}}$ and $\bigcup\limits^\bullet$ denotes a disjoint union.
According to the decomposition of the segment $D_{\kappa,\omega}(N)$ we introduce new Neumann boundary conditions,
which decreases  the operator. More precisely, we have in the sense of quadratic forms
\[
\mathcal{H}_{\kappa,\omega}(N) \geqslant \bigoplus_{j=0}^{J-1} \mathcal{H}_{\kappa,\omega}(K,j).
\]
In particular,
\begin{equation}
\label{e:Neumann-bracketing} \lambda_{\kappa,\omega}(N) \geqslant
\min_{j=0}^{J-1} \lambda_{\kappa,\omega}(K,j).
\end{equation}
The above considerations can be turned into a probabilistic
estimate on the position of the lowest eigenvalue.
Similar
ideas have been used e.g.~in \cite{MartinelliH-84,KirschSS-98a}
to obtain an initial scale estimate.
First note that by
(\ref{e:Neumann-bracketing}) we have the inclusion
\begin{align*}
\left\{\omega \in \Omega\mid \lambda_{\kappa,\omega}(N)-1\leqslant N^{-\frac{1}{2}}\right \}
&\subset
\bigcup_{j=0}^{J-1} \left\{\omega \in \Omega\mid \lambda_{\kappa,\omega}(K,j) -1\leqslant K^{-\frac{\gamma}{2}}\right \}
\end{align*}
Since the random variables $\omega_k, k \in \ZZ$ are independent and identically distributed,
\begin{align*}
\sum_{j=0}^{J-1} \PP\left(\omega \in \Omega\mid \lambda_{\kappa,\omega}(K,j) -1\leqslant K^{-\frac{\gamma}{2}}\right )
\leqslant
N^{1-\frac{1}{\gamma}} \ \PP\left(\omega \in \Omega\mid \lambda_{\kappa,\omega}(K) -1\leqslant K^{-\frac{\gamma}{2}}\right).
\end{align*}
By Corollary  \ref{th:lower-bound} and $\|\omega\|_{K,1} \leqslant \sqrt{K} \|\omega\|_{K,2}$
the following inclusions hold for all $\kappa \leqslant c_3 K^{-\frac{15}{2}}$:
\begin{align*}
 \left\{\omega \mid \lambda_{\kappa,\omega}(K) -1\leqslant  K^{-\frac{\gamma}{2}}\right\}
& \subset
 \left\{\omega \mid c_2 \frac{\kappa^2\|\omega\|_{K,2}^2}{K^3}\leqslant K^{-\frac{\gamma}{2}}\right\}
\\ =
  \left\{\omega \mid \|\omega\|_{K,2} \leqslant \frac{1}{\kappa\sqrt{c_2}} \, K^{\frac{3-\gamma/2}{2}}\right\}
& \subset
  \left\{\omega \mid \frac{\|\omega\|_{K,1} }{K}\leqslant \frac{1}{\kappa\sqrt{c_2}} \, K^{1-\frac{\gamma}{4}}\right\}.
\end{align*}
Denote by $\EE\{\omega_k\}$  the expectation value of (any) $\omega_k$,
and choose now $\kappa$ such that
\begin{equation}\label{e:kappa-conditions}
 \frac{K^{1-\frac{\gamma}{4}}}{\kappa\sqrt{c_2}} \leqslant \frac{\EE\{\omega_k\}}{2}
\quad \text{ i.e. } \quad
 \frac{2 K^{1-\frac{\gamma}{4}}}{\EE\{\omega_k\}\sqrt{c_2}} \leqslant \kappa.
\end{equation}
The upper and the lower bound for $\kappa$ can be reconciled if
$\gamma > 34$ and
\[
K \geqslant K_1 := \left( \frac{2}{\EE\{\omega_k\} c_3 \sqrt{c_2}}\right)^{\frac{2}{\gamma -34}}.
\]
The last inequality is equivalent to
$N\geqslant N_1:= ( \EE\{\omega_k\} c_3 \sqrt{c_2}/2)^{\frac{-2\gamma}{\gamma-34}}$.
For $\kappa$ satisfying (\ref{e:kappa-conditions}) we
are able to apply the large deviations principle of Lemma \ref{l:LDP} and thus obtain
\begin{align*}
\PP\left(\omega \mid \frac{\|\omega\|_1 }{K}\leqslant \frac{1}{\kappa\sqrt{c_2}} \, K^{1-\frac{\gamma}{4}}\right )
\leqslant
\PP\left(\omega \mid \frac{\|\omega\|_1 }{K}\leqslant \frac{\EE\{\omega_k\}}{2} \right)
\leqslant \E^{-c_4K}
\end{align*}
for $K\geqslant K_1$. Consequently we have for $N\geqslant N_1$
\[
\PP\left(\omega \in \Omega\mid \lambda_{\kappa,\omega}(N)-1\leqslant N^{-\frac{1}{2}}\right )
\leqslant
N^{1-\frac{1}{\gamma}} \ \E^{-c_4N^{1/\gamma}}.
\]
\end{proof}

\begin{proof}[Proof of Corollary  \ref{c:Combes-Thomas} ]
Since $ \lambda_{\kappa,\omega}(N)$ is the element of $\sigma(\mathcal{H}_{\kappa,\omega}(N))$ closest to one,
\begin{align*}
\Omega'
&=  \left\{\omega \in \Omega\mid \lambda_{\kappa,\omega}(N)-1> N^{-\frac{1}{2}}\right \}
=  \left\{\omega \in \Omega\mid \dist(1, \sigma(\mathcal{H}_{\kappa,\omega}(N) )> N^{-\frac{1}{2}} \right \}
\\
&=  \left\{\omega \in \Omega\mid \forall \ \lambda\in [1,1+1/(2\sqrt N)]  \colon
\dist(\lambda, \sigma(\mathcal{H}_{\kappa,\omega}(N) )> 1/(2\sqrt N) \right \}.
\end{align*}

For $\omega \in \Omega'$ we want to use a Combes-Thomas argument to obtain from the estimate on the
distance of the relevant energies to the spectrum a decay estimate for the Green's function.
Such estimates are rather standard for Schr\"odinger-type operators.  However, we did not find a
specific formulation of this result in the literature which fits exactly our situation, so we provide a direct proof,
for the convenience of the reader.
\medskip

In the following the parameters $\omega \in \Omega'$, $\kappa >0$,
and $N \in\NN$ will be kept fixed. For this reason we shall  
abbreviate in the subsequent calculations
$\mathcal{H}_{\kappa,\omega}(N)$ simply by $\mathcal{H}$.
Define the function $J \colon [0,l]\to [0, \infty)$ by setting
$J(t)=t$ for $ t \in (1,L-1)$,
$J(t)=3t^2-3
t^3+t^4$ for $ t \in [0,1]$,
and choosing the values on the segment $[L-1,L]$ in such a way that the graph of the function becomes point symmetric w.r.t. $(L/2,L/2)$.
Note that $J$ satisfies Neumann b.c.
at $t=0$ and $t=L$, that it is twice differentiable, and that on the segment  $ t \in [0,1]$
\[
J'(t)=  6t - 9 t^2 + 4 t^3 \geqslant 0, \quad
J''(t)= 6  -18 t   +12 t^2.
\]
It thus follows that
\begin{equation*}
\|J'\|_\infty \leqslant \frac{5}{4},
\quad
\|J''\|_\infty \leqslant 6,
\end{equation*}
and that $J$ is monotonously increasing.
For $a\in(0,1)$ we define the multiplication operator
\[
\mathcal{T}_a \colon L_2(D_{\kappa,\omega}(N))  \to L_2(D_{\kappa,\omega}(N)), \quad
(\mathcal{T}_a f)(x):= \E^{aJ(x_1)} f(x)
\]
and  another operator
\begin{equation*}
\mathcal{P}_a:=-2 a J' \frac{\partial}{\partial x_1}-a^2 J'-a J'',
\end{equation*}
which will turn out to be an `effective perturbation'.
A direct calculation shows that
\[
\mathcal{T}_{-a} \mathcal{T}_{a}=\I, \quad
\mathcal{T}_{-a} \mathcal{H} \mathcal{T}_{a}=\mathcal{H} + \mathcal{P}_a, \quad
\mathcal{T}_{-a} (\mathcal{H}- \lambda)^{-1} \mathcal{T}_{a}=(\mathcal{H} + \mathcal{P}_a - \lambda)^{-1},
\]
provided that $\lambda$ is in the intersection of the resolvent sets of the two operators $\mathcal{H}$
and $\mathcal{H}+\mathcal{P}_a$.
We 
 shall identify  a range of values for $a$ such that the last condition holds and
we get even an explicit bound on the norm of $(\mathcal{H} + \mathcal{P}_a - \lambda)^{-1}$.
In these considerations we 
shall keep in mind that $\lambda$ is close to, but larger that one, in particular $\lambda \in [1,2]$.
We denote by $\delta=\dist(\sigma(\mathcal{H}),\lambda)$ the distance of the spectrum of $\mathcal{H}$ to $\lambda$.
Let us first estimate
\begin{equation}
\begin{aligned}
\|\mathcal{P}_a (\mathcal{H}- \lambda)^{-1} \|
& \leqslant
a\bigg[ (a \|J'\|_\infty^2+\|J''\|_\infty) \|(\mathcal{H}- \lambda)^{-1} \|
+ 2 \|J'\|_\infty  \Big\|  \frac{\partial}{\partial x_1} (\mathcal{H}- \lambda)^{-1} \Big\| \bigg]
\\
& \leqslant
a\Bigg[ \frac{\frac{25}{16}a +  6}{\delta} +
 \frac{5}{2}\sqrt{\frac{\lambda}{\delta^2}+ \frac{1}{\delta}} \Bigg].
\end{aligned}\label{4.1}
\end{equation}
Here we have employed that
\begin{equation*}
\Big\|  \frac{\partial}{\partial x_1} (\mathcal{H}- \lambda)^{-1} \Big\| \leqslant \sqrt{\frac{\l}{\d^2}+\frac{1}{\d}},
\end{equation*}
which follows directly from the obvious relations
\begin{equation*}
\|\nabla u\|_{L_2(D_{\kappa,\om})}^2-\l \|u\|_{L_2(D_{\kappa,\om})}^2 = (f,u)_{L_2(D_{\kappa,\om})},\quad \|u\|_{L_2(D_{\kappa,\om})}\leqslant \frac{1}{\d}\|f\|_{L_2(D_{\kappa,\om})},
\end{equation*}
where $u:=(\mathcal{H}-\l)^{-1}f$, $f\in L_2(D_{\kappa,\om})$.

The
right hand side in (\ref{4.1}) is bounded by
\[
a \Bigg[\frac{121}{16\d}+\frac{5}{2}
 \sqrt{\frac{\lambda+\delta}{\delta^2}} \Bigg]
\leqslant 12 \frac{a}{\delta},
\]
since $\delta, a \in [0,1]$ and $\lambda \in [1,2]$.
Now choose $a = \frac{\delta}{24} \leqslant \frac{1}{2} \frac{\delta}{12}$. Then $\|\mathcal{P}_a (\mathcal{H}- \lambda)^{-1} \|
\leqslant 1/2$ and thus the norm of the Neumann series satisfies
\[
\|(\mathcal{H} + \mathcal{P}_a - \lambda)^{-1}\| \leqslant \frac{\|(\mathcal{H} - \lambda)^{-1}\| }{1-\frac 1 2 } = \frac{2}{\delta}.
\]
Now  choose $\alpha, \beta \geqslant2$ and set
\begin{align*}
A  &:=\{x \in \RR^2 \mid 0\leqslant  x_1\leqslant  \alpha  , \kappa G(x_1,\omega)<x_2<\kappa G(x_1,\omega)+\pi\} \subset D_{\kappa, \omega} (N),
\\
B  &:=\{x \in \RR^2 \mid L-\beta\leqslant  x_1\leqslant  L  , \kappa G(x_1,\omega)<x_2<\kappa G(x_1,\omega)+\pi\} \subset D_{\kappa, \omega} (N).
\end{align*}
For any normalised vectors $\phi,\psi \in L_2(D_{\kappa,\omega}(N))$ we have
\begin{align*}
\left \vert \left \la \vert \psi\vert \chi_A, \mathcal{T}_{-a} (\mathcal{H}-\lambda)^{-1}  \mathcal{T}_{a} \chi_B \vert \phi\vert \right \ra \right \vert
&\leqslant
 \|\mathcal{T}_{-a} (\mathcal{H}-\lambda)^{-1}  \mathcal{T}_{a}\|
=
\|(\mathcal{H} + \mathcal{P}_a - \lambda)^{-1}\|              \leqslant  \frac{2}{\delta}.
\end{align*}
Due to the monotonicity of $J$ and the positivity of the integral kernel of $(\mathcal{H}-\lambda)^{-1}$
 we are able to estimate
\begin{align*}
\left \vert \left \la \vert \psi\vert \chi_A, \mathcal{T}_{-a} (\mathcal{H}-\lambda)^{-1}  \mathcal{T}_{a} \chi_B \vert \phi\vert \right \ra \right \vert
&\geqslant
\E^{a(J(L-\beta)- J(\alpha))}
\left \vert \left \la  \psi  \chi_A,  (\mathcal{H}-\lambda)^{-1}  \chi_B  \phi \right \ra \right \vert.
\end{align*}
Note that by the choice of the values of $\alpha, \beta$  and the function $J$ we have $J(L-\beta)- J(\alpha) = L-\beta- \alpha
=\dist (A,B)$. Bringing the exponential term on the other side we obtain
\[
\left \vert \left \la  \psi  \chi_A,  (\mathcal{H}-\lambda)^{-1}   \chi_B  \phi \right \ra \right \vert
\leqslant
\frac{2}{\delta}  \E^{-a\dist(A,B)}
= \frac{2}{\delta}  \E^{-\frac{\dist(A,B)
\delta}{24}}
\]
Now fix  $\omega \in \Omega'=  \left\{\omega \in \Omega\mid \lambda_{\kappa,\omega}(N)-1> N^{-\frac{1}{2}}\right \}$
and $\lambda \in [1, 1+\frac{1}{2\sqrt N}]$.
Then $\delta \geqslant  \frac{1}{2\sqrt N}$ and thus
\[
\left \vert \left \la  \psi  \chi_A,  (\mathcal{H}-\lambda)^{-1}    \chi_B  \phi \right \ra \right \vert
\leqslant
\frac{2}{\delta}   \E^{-\frac{\dist(A,B)
 \delta}{
 24}}
\leqslant
2\sqrt N
\E^{-\frac{\dist(A,B)}{48\sqrt N}}.
\]
Since we have  by the estimate (\ref{eq:initial}) the bound $\PP(\Omega')\geqslant
1-N^{1-\frac{1}{\gamma}} \ \E^{-c_4N^{1/\gamma}}$ for $N \geqslant N_1$,  we conclude that
\begin{multline}
 \PP\left(\omega \in \Omega\mid \forall  \ \lambda\in [1,1+1/(2\sqrt N)]  \colon
\left \|\chi_A (\mathcal{H}_{\kappa,\omega}(N)-\lambda)^{-1}   \chi_B \right \|
\leqslant
2\sqrt N\E^{-\frac{\dist(A,B)}{48\sqrt N}}
\right)
\\
\geqslant
1-N^{1-\frac{1}{\gamma}} \ \E^{-c_4N^{1/\gamma}}.
\end{multline}
\end{proof}

\section{Deterministic lower bounds}
\label{s:deterministic}
In this section only finite segments of the waveguide will be relevant to us.
Recall that due to the fact that the support of the measure $\mu$
is contained in the interval $[0,1]$, for almost all $\omega \in \Omega$ the bound
$\|\omega\|_\infty \leqslant 1$ holds. In this section an $\ell^2$-normalisation will be better suited,
and in fact possible since only finite waveguide segments, and thus only finite families of
random variables $\omega_1, \dots, \omega_N$ are involved.

We fix $N \in \NN$ and set for all $j \in \{1,\dots,N\}$
\[
\rho_j :=\kappa \, \omega_j,
\theta_j:=\rho_j/\e,
\text{  where }
\e:=\left(\sum_{j=1}^N \rho_j^2\right)^{\frac{1}{2}}= \kappa\|\omega\|_{N,2}
\]
Observe that the vector $\theta=\{\theta_i\}_{i=1}^N $ is normalized in the sense
\begin{equation}\label{1.0a}
\theta_1^2+\ldots+\theta_N^2=1.
\end{equation}
Let $l\geqslant 1$, $g\colon \RR \to \RR$, $\widetilde{g}\in \RR$, $G\colon\RR \times \Omega \to \RR$ be as in Sections \ref{s:model} and \ref{s:probab-shifted}.
We set $L:=Nl$ and consider in the following the restriction of $G$
to the set $[0,L]\times \big(\times_{j=1}^N [0,\infty)\big)$.
The restricted function will be again denoted by $G$, and the finite vector $(\omega_1, \dots,\omega_N)$
will be denoted by $\omega$, by slight abuse of notation.
In terms of the new parameters we have
\[
\kappa G(x_1, \omega)= G(x_1, \kappa \omega)= G(x_1, \rho)= G(x_1, \e \theta)
= \e G(x_1, \theta)
=\e \sum_{j=0}^N\theta_{j} g(x_1-jl).
\]
For the subset of $\RR^2$ forming the finite segment of the waveguide we 
shall use in this section the notation
\[
\Pi_\rho:=\{x\in \RR^2 \mid  0<x_1<L, \e G(x_1,\theta)<x_2<\e G(x_1,\theta)+\pi\}
\]
Let us point out  that $\Pi_\rho= D_{\|\rho\|_\infty,\rho/\|\rho\|_\infty}(N)$ in the notation of Section \ref{s:model}
and that $\Pi_0= \{x\in \RR^2 \mid  0<x_1<L, 0<x_2<\pi\}$.
Similarly as before
\begin{equation*}
\G_\rho:=\{x\in\RR^2 \mid x_1\in(0,L), x_2=\e G(x_1,\theta)\}\cup\{x: x_1\in(0,L),
x_2=\e G(x_1,\theta)+\pi\}.
\end{equation*}
indicates the upper and lower parts of the boundary of $\Pi_\rho$ and the remainder of the boundary $\p\Pi_\rho\setminus\overline{\G}_\rho$ is denoted by $\g_\rho$. By $\mathcal{H}_\rho$ we denote the negative Laplacian in $\Pi_\rho$ with Dirichlet boundary conditions on $\G_\rho$ and Neumann ones on $\g_\rho$ and by $\lambda(\rho)$ the lowest eigenvalue of $\mathcal{H}_\rho$.

In what follows $\la\cdot,\cdot\ra$ and $\|\cdot\|$  denote  the scalar product
in $L_2(\Pi_0)$ and the associated norm.
We repeat here the statement of Theorem \ref{th:deterministic-lower-bound-1} in the new parametrization.

\begin{theorem}
For
\begin{equation}\label{1.0c}
\e  \leqslant   \frac{3}{5000}  \, \|\widetilde{g}\|_{L_2(0,l)}^2  \, \frac{1}{ L^7}
\end{equation}
the estimate
\begin{equation*}
\lambda(\rho)-1\geqslant
\frac{3}{2}\, \|\widetilde{g}\|_{L_2(0,l)}^2 \, \frac{\e^2}{L^3}
\end{equation*}
holds true.
\end{theorem}

The rest of the section is devoted to the proof of the theorem.
\bigskip 

First we transform the Laplacian on the 
wiggled strip segments into a differential operator with variable coefficients on a rectangle.
It is then possible to treat the latter operator as a perturbation of the pure Laplacian on the rectangle.
We introduce the coordinates $\xi:=(\xi_1,\xi_2)$, $\xi_1:=x_1$,
$\xi_2:=x_2-\e G(x_1)$. The mapping $u(x)\mapsto u(\xi_1,\xi_2+\e
G(\xi_1))$ is a unitary operator from $L_2(\Pi_\rho)$ to
$L_2(\Pi_0)$. Let $\psi_\rho$ be the normalized eigenfunction
associated with $\lambda(\rho)$. It is the unique solution of the boundary value problem
\begin{gather}
(-\D_\xi-\e
\mathcal{Q}_\rho)\psi_\rho=\lambda(\rho)\psi_\rho,\quad
\xi\in\Pi_0,
\\
\psi_0=0,\quad x\in\G_0,\qquad \frac{\p\psi_0}{\p\xi_1}=0,\quad
\xi\in\g_0,
\end{gather}
where
\begin{equation*}
\mathcal{Q}_\rho:=-2G'\frac{\p^2}{\p\xi_1\p\xi_2}+\e
\big(G'\big)^2\frac{\p^2}{\p\xi_2^2}-G''\frac{\p}{\p\xi_2}.
\end{equation*}
Hereafter the prime denotes the derivative w.r.t. $\xi_1$.
 \medskip

For the perturbation theoretic estimates we are aiming for  we need some control of the resolvent set.
The first two eigenvalues of $\mathcal{H}_0$ are $1$ and
$1+\pi^2/L^2$.  The eigenfunction corresponding to the eigenvalue one is
\[
\psi_0\colon \Pi_0\to \Pi_0, \quad  \psi_0(\xi):=\sqrt{\frac{2}{\pi L}}\sin\xi_2.
\]
It turns out that it is more convenient to work with the modified resolvent
\begin{equation*} 
\mathcal{R}_0(\lambda)
:=(\mathcal{H}_0-\l)^{-1}
-\frac{\la\cdot, \psi_0\ra}{1-\l}\psi_0.
\end{equation*}
which is a bounded self-adjoint operator for any
\begin{equation}
\label{1.5}
\lambda\in B_{\pi^2/(2L^2)}(1)\subset \CC.
\end{equation}
The operator-valued function $B_{\pi^2/(2L^2)}(1)\ni \lambda \mapsto \mathcal{R}_0(\lambda)$ is holomorphic, cf.{} \cite[Ch. V, Sec. 3.5]{Kato-66}

For an arbitrary $f\in L_2(\Pi_0)$ set $\widehat{f}:=f-\la f,\psi_0\ra \psi_0$. Then the function $u:=\mathcal{R}_0(\lambda)f= (\mathcal{H}_0-\l)^{-1} \widehat{f}$ solves the equation
\begin{equation}\label{1.6}
(\mathcal{H}_0-\l)u
= (\mathcal{H}_0-\l)\, (\mathcal{H}_0-\l)^{-1} f + (\mathcal{H}_0-\l) \frac{\la f , \psi_0\ra }{\l -1}\psi_0=
\widehat{f},
\end{equation}
 and is orthogonal to $\psi_0$ in $L_2(\Pi_0)$.

\begin{lemma}\label{lm1.1}
For $f \in L_2(\Pi_0)$, $u:=\mathcal{R}_0(\lambda)f$ and
$\l\in B_{\pi^2/(2L^2)}(1)\subset \mathbb{C}$ we have $
u\in W^2_2(\Pi_0)$ and
\begin{equation}
\begin{aligned}
&\|u\|  \leqslant   \frac{2L^2}{\pi^2} \|f\|,&&\|\nabla u\|  \leqslant
\frac{7 L^2}{\pi^2}\|f\|,
\\
& \Big\|\nabla\frac{\p u}{\p\xi_1}\Big\|  \leqslant   \frac{25 L^2}{\pi^2}\|f\|, && \Big\|\frac{\p^2
u}{\p\xi_2^2}\Big\|  \leqslant  \frac{47 L^2}{\pi^2}\|f\|\label{1.7}.
\end{aligned}
\end{equation}
\end{lemma}

\begin{proof}
The vector $u= (\mathcal{H}_0-\l)^{-1} \widehat{f}$
is in the range of  $(\mathcal{H}_0-\l)^{-1}$ and thus in the Sobolev space $W^2_2(\Pi_0)$.
The first inequality follows from \cite[Ch.\,V]{Kato-66} and the fact that our choice of $\lambda$ is separated by at least $\pi^2/(2L^2)$ from the next spectral value.

Let us prove the second estimate.  We begin with the obvious identity
\begin{equation}\label{1.8}
\|\widehat{f}\|^2=\|f\|^2-|\la f,\psi_0\ra^2|  \leqslant   \|f\|^2.
\end{equation}
We multiply (\ref{1.6}) with $u$ and obtain $\|\nabla u\|^2=\l\|u\|^2+\la\widehat{f},u\ra$.
This identity, (\ref{1.5}) and the first estimate in (\ref{1.7}) yield
\begin{align}
&
\begin{aligned}
\|\nabla u\|^2&
\leqslant |\lambda| \|u\|^2+\|f\|\,\|u\|
\leqslant \frac{4L^4}{\pi^4}|\lambda| \|f\|^2+\frac{2L^2}{\pi^2}\|f\|^2
\\
&\leqslant \frac{4L^4}{\pi^4}\left(1+ \frac{\pi^2}{2L^2}+\frac{\pi^2}{2L^2}\right)\|f\|^2
 \leqslant
\frac{4L^4}{\pi^4}\left(1+\pi^2\right)\|f\|^2
\end{aligned}\label{1.9}
\end{align}
that proves the desired estimate.

Since $u\in \H^1(\Pi_0)$, the function $v:=\frac{\p u}{\p\xi_1}$ is a generalized solution to the boundary value problem
\begin{equation*}
-\D_\xi v=\l v+\frac{\p\widehat{f}}{\p\xi_1},\quad
\xi\in\Pi_0, \qquad
 v=0,\quad
 \xi\in\p\Pi_0
\end{equation*}
which is obtained by differentiating equation (\ref{1.6}).
We multiply the last equation by $v$, integrate by parts,
and obtain
\[
\|\nabla v\|^2=\l\|v\|^2-\Big\la \widehat{f}, \frac{\p v}{\p\xi_1}\Big\ra. ,
\]
We employ (\ref{1.5}), (\ref{1.8}), and Young's inequality  to estimate
\[
 \|\nabla v\|^2
  \leqslant   \left(1+\frac{\pi^2}{2 L^2}\right)\|v\|^2+\frac{1}{2}\|f\|^2+\frac{1}{2}\Big\| \frac{\p v}{\p\xi_1}\Big\|^2
  \leqslant   \left(1+\frac{\pi^2}{2 L^2}\right)\|\nabla u \|^2+\frac{1}{2}\|f\|^2+\frac{1}{2} \|\nabla  v \|^2,
\]
and hence by (\ref{1.9})
\begin{equation}
\label{1.10}
\begin{aligned}
\Big\|\nabla\frac{\p u}{\p\xi_1}\Big\|^2&=\|\nabla v\|^2
\leqslant
\left(2+\frac{\pi^2}{L^2}\right)\|\nabla u\|^2+\|f\|^2
\\
&\leqslant \left(4(\pi^2+1)\left(2+\frac{\pi^2}{L^2}\right)+\frac{\pi^4}{L^4}
\right) \frac{L^4}{\pi^4}\|f\|^2
\\
&\leqslant \big(4(\pi^2+1) (2+\pi^2 )+ \pi^4
\big) \frac{L^4}{\pi^4}\|f\|^2\leqslant \frac{614L^4}{\pi^4}\|f\|^2.
\end{aligned}
\end{equation}
This implies the penultimate inequality in (\ref{1.7}).

We rewrite (\ref{1.6}) as
\begin{equation*}
-\frac{\p^2 u}{\p\xi_2^2}= \frac{\p^2 u}{\p\xi_1^2}+\l
u+\widehat{f},
\end{equation*}
and see that the  estimates (\ref{1.5}), (\ref{1.8}), (\ref{1.10}),  and the first estimate
in (\ref{1.7}) yield
\begin{align*}
\Big\|\frac{\p^2 u}{\p\xi_2^2}\Big\| &  \leqslant   \Big\|\frac{\p^2 u}{\p\xi_1^2}\Big\|+|\lambda|\, \|u\|+\|f\|
\\
&\leqslant  \left(\sqrt{614}+2\left(1+\frac{\pi^2}{2 L^2}\right)+\pi^2\right) \frac{L^2}{\pi^2} \|f\|
\\
&\leqslant  \left(\sqrt{614}+2 (1+ \pi^2)+\pi^2\right) \frac{L^2}{\pi^2} \|f\|\leqslant \frac{47 L^2}{\pi^2}\|f\|,
\end{align*}
which proves the last claim in (\ref{1.7}).
\end{proof}
The inequalities in (\ref{lm1.1}) can be understood as bounds on the norms of certain operators products:
\begin{equation}
\label{eq:operator-norms}
\begin{aligned}
&\|\mathcal{R}_0(\lambda)  \|\leqslant   \frac{2L^2}{\pi^2},
&&\|\nabla \mathcal{R}_0(\lambda)\|  \leqslant
\frac{7 L^2}{\pi^2},
\\
& \Big\|\nabla\frac{\p }{\p\xi_1}\mathcal{R}_0(\lambda)\Big\|  \leqslant   \frac{25 L^2}{\pi^2},
&& \Big\|\frac{\p^2}{\p\xi_2^2}\mathcal{R}_0(\lambda)\Big\|  \leqslant  \frac{47 L^2}{\pi^2}.
\end{aligned}
\end{equation}

By the normalization (\ref{1.0b}) we have
\begin{equation}\label{1.10a}
\|\widetilde{g}\|_{L_2(0,l)}^2=\|g\|_{L_2(0,l)}^2-l^{-1}\Big(\int\limits_0^l
g(t)\di t\Big)^2  \leqslant   \|g\|_{L_2(0,l)}^2  \leqslant   l.
\end{equation}
This implies together with inequality (\ref{1.0c})
\begin{equation}\label{1.10b}
\e
\leqslant   \frac{3l^2}{5000 L^7}
\leqslant   \frac{3}{5000}   \frac{1}{ L^5\, N^2}
\leqslant  \frac{3}{5000}.
\end{equation}

We use (\ref{1.0a}) and (\ref{1.0b}) and Cauchy-Schwartz inequality to establish
\begin{equation*}
\|\mathcal{Q}_\rho u \|
\leqslant 2\Big\|\frac{\p^2 u }{\p\xi_1\p\xi_2}\Big\|+ \Big\|\frac{\p^2 u }{\p\xi_2^2}\Big\|+ \Big\|\frac{\p u }{\p\xi_2}\Big\|
\end{equation*}
for any $ u \in W_2^2(\Pi_0)$. For  $\lambda$ satisfying (\ref{1.5}), and $\e$ as in  (\ref{1.10b}),
Lemma~\ref{lm1.1} and the last bound
imply
\begin{equation}\label{1.12}
\|\mathcal{Q}_\rho \mathcal{R}_0(\lambda)\|
\leqslant   \frac{
 104 L^2}{\pi^2},\quad \e\|\mathcal{Q}_\rho \mathcal{R}_0(\lambda)\|
\leqslant   \frac{
8}{125\pi^2}<1.
\end{equation}
Here the norm is understood as the norm of an operator in $L_2(\Pi_0)$.

Let us show that $\lambda(\rho)$ satisfies (\ref{1.5}).  Indeed,
\begin{equation}\label{1.3}
\lambda(\rho)=\Big\|\frac{\p \psi_\rho}{\p\xi_1}-\e G'\frac{\p
\psi_\rho}{\p \xi_2}\Big\|^2+\Big\|\frac{\p
\psi_\rho}{\p\xi_2}\Big\|^2\geqslant \Big\|\frac{\p
\psi_\rho}{\p\xi_2}\Big\|^2\geqslant 1.
\end{equation}
By the minimax principle and (\ref{1.0a}), (\ref{1.0b}) the
eigenvalue $\lambda(\rho)$ can be estimated from above as follows
\begin{equation}\label{1.3a}
\lambda(\rho)  \leqslant   \Big\|\frac{\p \psi_0}{\p\xi_1}-\e G'\frac{\p
\psi_0}{\p \xi_2}\Big\|^2+\Big\|\frac{\p
\psi_0}{\p\xi_2}\Big\|^2  \leqslant   1+\e^2.
\end{equation}
Two last estimates and (\ref{1.10b}) imply (\ref{1.5}) for
$\lambda(\rho)$.
\bigskip

Next we 
shall do some perturbation theory for linear operators to be
able to identify the dominating contributions. Using the proved
fact, (\ref{1.12}), and proceeding completely as in
\cite{Gadylshin-02}, \cite[Sec. 4]{Borisov-06}, one can show that
for the considered values of $\e$ the eigenvalue $\lambda(\rho)$
solves the equation
\begin{equation}\label{1.13}
\lambda(\rho)-1=-\e\big\la (\I-\e \mathcal{Q}_\rho
\mathcal{R}_0(\lambda(\rho)))^{-1}\mathcal{Q}_\rho\psi_0,\psi_0\big\ra,
\end{equation}
where $\I$ is the identity mapping. A direct calculation shows
\begin{equation*}
(\I-\e \mathcal{Q}_\rho \mathcal{R}_0(\lambda(\rho)))^{-1}=\I+\e
\mathcal{Q}_\rho \mathcal{R}_0(\lambda(\rho))+\e^2 (\I-\e
\mathcal{Q}_\rho \mathcal{R}_0(\lambda(\rho)))^{-1}
\big(\mathcal{Q}_\rho \mathcal{R}_0(\lambda(\rho))\big)^2.
\end{equation*}
We substitute this relation into (\ref{1.13}),
\begin{equation}\label{1.14}
\begin{aligned}
\lambda(\rho)-1=&-\e \la \mathcal{Q}_\rho\psi_0,\psi_0\ra-\e^2
\la \mathcal{Q}_\rho \mathcal{R}_0(\l(\rho)) \mathcal{Q}_\rho
\psi_0,\psi_0\ra
\\
&-\e^3 \big\la (\I-\e \mathcal{Q}_\rho
\mathcal{R}_0(\lambda(\rho)))^{-1} \big(\mathcal{Q}_\rho
\mathcal{R}_0(\lambda(\rho))\big)^2 \mathcal{Q}_\rho
\psi_0,\psi_0\big\ra.
\end{aligned}
\end{equation}
We first consider the coefficient of $\e^2$. For the sake of brevity set
\begin{equation*}
 v :=\mathcal{R}_0(\lambda(\rho)) \mathcal{Q}_\rho\psi_0.
\end{equation*}
Integrating by parts, we obtain
\begin{align*}
\la \mathcal{Q}_\rho \mathcal{R}_0(\lambda(\rho)) \mathcal{Q}_\rho\psi_0,\psi_0\ra
=& \la \mathcal{Q}_\rho  v ,\psi_0\ra
=-2\Big\la G' \frac{\p^2  v }{\p\xi_1\p\xi_2},\psi_0\Big\ra-\Big\la G'' \frac{\p  v }{\p\xi_2},\psi_0\Big \ra + \e \Big\la (G')^2 \frac{\p^2  v }{\p\xi_2^2},\psi_0\Big\ra
\\
=&-2\Big\la G''  v ,\frac{\p\psi_0}{\p\xi_2}\Big\ra+\Big\la G''  v ,\frac{\p\psi_0}{\p\xi_2}\Big \ra + \e \Big\la (G')^2  v ,\frac{\p^2 \psi_0}{\p\xi_2^2}\Big\ra
\\
=&-\Big\la G''  v ,\frac{\p\psi_0}{\p\xi_2}\Big\ra- \e \Big\la (G')^2  v ,\psi_0\Big\ra
\\
=& -\Big\la G'' \mathcal{R}_0(\lambda(\rho)) \mathcal{Q}_\rho\psi_0,
\frac{\p\psi_0}{\p\xi_2}\Big\ra -\e \la (G')^2\mathcal{R}_0(\lambda(\rho))
\mathcal{Q}_\rho\psi_0,\psi_0\ra
\end{align*}
 which equals
\begin{align*}
 \Big\la G'' \mathcal{R}_0(\lambda(\rho)) G''\frac{\p\psi_0}{\p\xi_2},
\frac{\p\psi_0}{\p\xi_2} \Big\ra + \e\Big\la G''
\mathcal{R}_0(\lambda(\rho)) (G')^2\psi_0, \frac{\p\psi_0}{\p\xi_2}
\Big\ra
 -\e \la (G')^2\mathcal{R}_0(\lambda(\rho))
\mathcal{Q}_\rho\psi_0,\psi_0\ra.
\end{align*}
Since $\mathcal{R}_0(\lambda)$ is  a (special) rank one perturbation of $ (\mathcal{H}_0-\l)^{-1}$, cf. (\ref{1.6}), it inherits the resolvent identity
\begin{equation*}
\mathcal{R}_0(\lambda)-\mathcal{R}_0(\mu)=(\l-\mu) \mathcal{R}_0(\mu)
\mathcal{R}_0(\lambda).
\end{equation*}
Now we substitute two last identities into (\ref{1.14}),
\begin{equation}\label{1.15}
\begin{aligned}
\lambda(\rho)-1=&-\e \la \mathcal{Q}_\rho\psi_0,\psi_0\ra-
\e^2\Big\la G'' \mathcal{R}_0(1) G''\frac{\p\psi_0}{\p\xi_2},
\frac{\p\psi_0}{\p\xi_2} \Big\ra
\\
&+\e^3 \la (G')^2\mathcal{R}_0(\lambda(\rho))
\mathcal{Q}_\rho\psi_0,\psi_0\ra
\\
& - \e^3\Big\la G'' \mathcal{R}_0(\lambda(\rho)) (G')^2\psi_0,
\frac{\p\psi_0}{\p\xi_2} \Big\ra
\\
&-\e^2(\lambda(\rho)-1) \Big\la G''
\mathcal{R}_0(1)\mathcal{R}_0(\lambda(\rho))
G''\frac{\p\psi_0}{\p\xi_2}, \frac{\p\psi_0}{\p\xi_2} \Big\ra
\\
&-\e^3 \big\la (\I-\e \mathcal{Q}_\rho
\mathcal{R}_0(\lambda(\rho)))^{-1} \big(\mathcal{Q}_\rho
\mathcal{R}_0(\lambda(\rho))\big)^2\mathcal{Q}_\rho\psi_0,\psi_0\big\ra.
\end{aligned}
\end{equation}

It will turn out that the two first terms on the right hand side of this
identity are positive and dominate all the remaining terms.
The first term is of the order $\e^2/L$. However, a cancelation occurs such that the
leading contribution turns out to be  of order $\e^2/L^3$.
A direct calculation shows that
\begin{equation}\label{1.16}
\la
\mathcal{Q}_\rho\psi_0,\psi_0\ra=-\frac{\e}{L}\|g'\|_{L_2(0,l)}^2.
\end{equation}
It is more complicated to evaluate the second term. We denote
\begin{equation*}
 w :=\mathcal{R}_0(1) G''\frac{\p\psi_0}{\p\xi_2},
\end{equation*}
Since this function solves the equation
\begin{equation*}
(\mathcal{H}_0-1)w = G''\frac{\p\psi_0}{\p\xi_2},
\end{equation*}
it is natural to consider the auxiliary function
\begin{equation}\label{1.17a}
v:= w +G\frac{\p\psi_0}{\p\xi_2}.
\end{equation}
Using the normalization (\ref{1.0a}), separation of variables and
\begin{equation*}
\| \frac{\partial \psi_0}{\partial \xi_2}\|^2_{L_2(0,\pi)}=1/L,
\end{equation*}
 we get
\begin{equation}\label{1.18}
\Big\la G'' \mathcal{R}_0(1) G''\frac{\p\psi_0}{\p\xi_2},
\frac{\p\psi_0}{\p\xi_2} \Big\ra= \Big\la G''  w ,
\frac{\p\psi_0}{\p\xi_2} \Big\ra=\frac{1}{L}\|g'\|_{L_2(0,l)}^2+
\Big\la v,G''\frac{\p\psi_0}{\p\xi_2} \Big\ra.
\end{equation}
The first term on the right side  leads to a cancelation of (\ref{1.16}).
It is clear that the function $v$ solves the boundary value
problem
\begin{gather*}
-\D v=v,\quad \xi\in\Pi_0,\qquad \frac{\p v}{\p\xi_1}=0,\quad
\xi\in\g_0,
\\
v=G\sqrt{\frac{2}{\pi L}},\quad \xi_1\in(0,L),\ \xi_2=0, \qquad
v=-G\sqrt{\frac{2}{\pi L}},\quad \xi_1\in(0,L),\ \xi_2=\pi.
\end{gather*}
We multiply the equation in this problem by
$G\frac{\p\psi_0}{\p\xi_2}$ and integrate twice by parts.
Due to the separation of variables a direct calculation shows that part of the boundary contributions vanish and the equation simplifies to
\begin{align}
\Big\la v,G\frac{\p\psi_0}{\p\xi_2}\Big\ra=&\sqrt{\frac{2}{\pi
L}}\int\limits_0^L G(\xi_1,\theta)\left(\frac{\p
v}{\p\xi_2}\Big|_{\xi_2=\pi}+\frac{\p
v}{\p\xi_2}\Big|_{\xi_2=0}\right)\di\xi_1\nonumber
\\
&-\Big\la v,G''\frac{\p\psi_0}{\p\xi_2}\Big\ra+\Big\la
v,G\frac{\p\psi_0}{\p\xi_2}\Big\ra. \nonumber 
\intertext{Thus}
\Big\la v,G''\frac{\p\psi_0}{\p\xi_2}\Big\ra=&
\sqrt{\frac{2}{\pi L}}\int\limits_0^L G(\xi_1,\theta)\left(\frac{\p
v}{\p\xi_2}\Big|_{\xi_2=\pi}+\frac{\p
v}{\p\xi_2}\Big|_{\xi_2=0}\right)\di\xi_1.\label{1.19}
\end{align}
We 
 shall now expand the functions in terms of the eigenvectors of the transversal Laplacian with Dirichlet b.c.
and the longitudinal Laplacian with Neumann b.c.
It follows from the equation
$(\mathcal{H}_0-\l)w=G''\frac{\p\psi_0}{\p\xi_2}$ that
\begin{equation*}
 w (\xi,\rho)=-8\pi\sqrt{\frac{2}{\pi L}}
\sum\limits_{n,m=1}^{\infty} \frac{n^2 m G_n(\rho)}{4m^2-1}
\frac{\cos \pi n L^{-1}\xi_1\sin 2m\xi_2}{(4m^2-1)L^2+\pi^2 n^2},
\end{equation*}
where the coefficients $G_n(\rho)$ are defined by the identity
\begin{equation*}
G(\xi_1,\theta)=\sum\limits_{n=0}^{\infty} G_n(\rho)\cos \pi n
L^{-1} \xi_1.
\end{equation*}
We substitute the expansion obtained and (\ref{1.17a}) into
(\ref{1.19}),
\begin{align*}
-\Big\la v,G''\frac{\p\psi_0}{\p\xi_2}  \Big\ra
&=32
\sum\limits_{n,m=1}^{\infty} \frac{n^2 m^2 G_n^2(\rho)}{4m^2-1}
\frac{1}{(4m^2-1)L^2+\pi^2 n^2}
\\
 &>\frac{32}{3} \sum\limits_{n=1}^{\infty} \frac{n^2
G_n^2(\rho)}{3 L^2+\pi^2 n^2}
>
\frac{32}{3(3+\pi^2)L^2}
 \sum\limits_{n=1}^{\infty}
G_n^2(\rho)>
\frac{64}{(9+3\pi^2)L^2}
\|\widetilde{g}\|_{L_2(0,l)}^2.
\end{align*}
Here we have also employed that due to (\ref{1.10a}) and the normalization (\ref{1.16})
\begin{align*}
\sum\limits_{n=1}^{\infty} \frac{L G_n^2(\rho)}{2}&=\|G-
G_0\|_{L_2(0,L)}^2=\|G\|_{L_2(0,L)}^2-L G_0^2
\\
&=\|g\|_{L_2(0,l)}^2 -\frac{1}{L}
\left(\sum\limits_{i=1}^{N}\theta_i\right)^2\left(\int\limits_0^l
g(t)\di t\right)^2
\\
&\geqslant \|g\|_{L_2(0,l)}^2 -\frac{1}{l} \left(\int\limits_0^l
g(t)\di t\right)^2=\|\widetilde{g}\|_{L_2(0,l)}^2.
\end{align*}

Now it follows from (\ref{1.16}), (\ref{1.18}) that
\begin{equation}\label{1.20}
-\e \la \mathcal{Q}_\rho\psi_0,\psi_0\ra
- \e^2\Big\la G'' \mathcal{R}_0(1) G''\frac{\p\psi_0}{\p\xi_2}, \frac{\p\psi_0}{\p\xi_2} \Big\ra
\geqslant \frac{32\e^2}{3\pi^2 L^3 }\|\widetilde{g}\|_{L_2(0,l)}^2.
\end{equation}
Using (\ref{1.0b}), (\ref{eq:operator-norms}), (\ref{1.3}),
(\ref{1.3a}), $\| \frac{\partial \psi_0}{\partial \xi_2}\|=1 $, and the inequality
\begin{equation*}
\|\mathcal{Q}_\rho\psi_0\|  \leqslant   1+\e  \leqslant   2,
\end{equation*}
we estimate the remaining terms in the right hand side of (\ref{1.15}) by noting that
\begin{align*}
&\e^3|\la (G')^2\mathcal{R}_0(\lambda(\rho))
\mathcal{Q}_\rho\psi_0,\psi_0\ra|  \leqslant
\frac{4\e^3L^2}{\pi^2}  \leqslant   \frac{4\e^3 L^4}{\pi^2},
\\
& \e^3\Big|\Big\la G'' \mathcal{R}_0(\lambda(\rho)) (G')^2\psi_0,
\frac{\p\psi_0}{\p\xi_2} \Big\ra\Big|  \leqslant
\frac{2\e^3L^2}{\pi^2}
\leqslant   \frac{2\e^3L^4}{\pi^2}.
\end{align*}
To estimate the penultimate  term (\ref{1.15}) observe that by
(\ref{1.3}), (\ref{1.3a}) and (\ref{1.10b})
\[
\Big|\lambda(\rho)-1) \Big| \leqslant \e^2 \leqslant \frac{3\e}{5000 L^5 N^2}
\leqslant \frac{3\e}{5000 L^5}
\]
since $N\geqslant 1$. The last estimate, (\ref{1.0b}), and the first inequality in
(\ref{eq:operator-norms})
yield
\begin{multline*}
\e^2\Big|(\lambda(\rho)-1) \Big\la G'' \mathcal{R}_0(1)\mathcal{R}_0(\lambda(\rho))
G''\frac{\p\psi_0}{\p\xi_2}, \frac{\p\psi_0}{\p\xi_2}\Big\ra\Big|
\leqslant
\frac{3\e^3}{5000 L^5} \Big( \frac{2L^2}{\pi^2}\Big)^2 \Big\| \frac{\partial \psi_0}{\partial \xi_2}\Big\|_{L_2(\Pi_0)}^2
\\
\leqslant
\frac{3\e^3}{1250\pi^4 L}  \leqslant
\frac{3\e^3 L^4}{1250\pi^4}.
\end{multline*}
We use (\ref{1.12}) and the estimate
\begin{equation*}
 \|\mathcal{Q}_\rho\psi_0\|  \leqslant   1+\e \leqslant   2
\end{equation*}
to arrive at
\begin{align*}
&\e^3\big|\big\la (\I-\e \mathcal{Q}_\rho
\mathcal{R}_0(\lambda(\rho)))^{-1} \big(\mathcal{Q}_\rho
\mathcal{R}_0(\lambda(\rho))\big)^2
\mathcal{Q}_\rho\psi_0,\psi_0\big\ra\big|
\leqslant
\frac{250\cdot 104^2 L^4}{125\pi^4-8\pi^2}\e^3.
\end{align*}
In view of (\ref{1.15}), (\ref{1.20}) it leads us to the
estimate
\begin{equation}\label{eq:ultimate}
\lambda(\rho)-1
\geqslant
\frac{64\e^2}{(9+3\pi^2)L^3} \|\widetilde{g}\|_{L_2(0,l)}^2- 225 L^4\e^3.
\end{equation}
We note that
\[
\frac{64}{9+3\pi^2} -3\frac{225}{5000} > \frac{3}{2}
\]
and use  assumption  (\ref{1.0c}) to bound (\ref{eq:ultimate})
from below by
\[
\frac{64\e^2}{(9+3\pi^2)L^3} \|\widetilde{g}\|_{L_2(0,l)}^2- 225 L^4\e^2 \frac{3\|\widetilde{g}\|_{L_2(0,l)}^2}{5000 L^7}
\geqslant \frac{3}{2} \|\widetilde{g}\|_{L_2(0,l)}^2 \frac{\e^2}{L^3},
\]
which completes the proof.

\section*{Acknowledgments}

D.B. was partially supported by RFBR,
by the grants of the President of Russia for young
scientists-doctors of sciences (MD-453.2010.1) and for Leading
Scientific Schools (NSh-6249.2010.1), and by the Federal Task Program ``Scientific and Scientific-Pedagogical Personnel of the Innovative Russia in 2009-2013'' (contract No. 02.740.11.0612).
I.V. was partially supported by the Deutsche Forschungsgemeinschaft
through the project
`Spectral properties of random Schr\"odinger operators and random operators on manifolds and graphs'
within the Emmy-Noether-Programme.

\def\cprime{$'$}\def\polhk#1{\setbox0=\hbox{#1}{\ooalign{\hidewidth
  \lower1.5ex\hbox{`}\hidewidth\crcr\unhbox0}}}

\end{document}